\documentclass{amsart}

\usepackage{amsmath,amssymb,amsthm}
\usepackage{graphicx,subcaption,float,caption}

\usepackage{graphicx,tikz}

\newtheorem*{thm}{Theorem}
\newtheorem*{proposition}{Proposition}

\newtheorem*{lem}{Lemma}

\theoremstyle{definition}

\theoremstyle{remark}

\begin{document}

\title[]{A Semicircle Law for derivatives\\ of random polynomials}


\author[]{Jeremy G. Hoskins}
\address{Program in Applied Mathematics, Yale University, New Haven, CT 06511, USA}
\email{jeremy.hoskins@yale.edu}

\author[]{Stefan Steinerberger}
\address{Department of Mathematics, Yale University, New Haven, CT 06511, USA}
\email{stefan.steinerberger@yale.edu}
\thanks{S.S. is supported by the NSF (DMS-1763179) and the Alfred P. Sloan Foundation.}

\begin{abstract} 
Let $x_1, \dots, x_n$ be $n$ independent and identically distributed random variables with mean zero, unit variance, and finite moments of all remaining orders. We study the random polynomial $p_n$ having roots at $x_1, \dots, x_n$.
We prove that for $\ell \in \mathbb{N}$ fixed as $n \rightarrow \infty$, the $(n-\ell)-$th derivative of $p_n^{}$ behaves like a Hermite polynomial: for $x$ in a compact interval,
 $${n^{\ell/2}} \frac{\ell!}{n!} \cdot p_n^{(n-\ell)}\left( \frac{x}{\sqrt{n}}\right) \rightarrow He_{\ell}(x + \gamma_n),$$
 where $He_{\ell}$ is the $\ell-$th probabilists' Hermite polynomial and $\gamma_n$ is a random variable converging to the standard $\mathcal{N}(0,1)$ Gaussian as $n \rightarrow \infty$. Thus, there is a universality phenomenon when differentiating a random polynomial many times: the remaining roots follow a Wigner semicircle distribution.
\end{abstract}

\maketitle

\section{Introduction}
\subsection{Introduction.} Let $p_n:\mathbb{R} \rightarrow \mathbb{R}$ be a polynomial of degree $n$ having
$n$ distinct real roots. Rolle's theorem implies that the derivative $p_n'$ has exactly
$n-1$ real roots. Moreover, between any two roots of $p_n$ there is exactly one root
of $p_n'$, the roots interlace. We are motivated by the following question \cite{steini}.
\begin{quote} \textbf{Open Problem.} Let $x_1, \dots, x_n$ be $n$ i.i.d. random variables and let $p_n$ be a polynomial of degree $n$ vanishing at these $n$ points.
What can be said about the distribution of roots of the $(t \cdot n)-$th derivative when $0 < t < 1$?
\end{quote}
While the question is interesting in and of itself, related problems also arise in the spectra of restrictions of symmetric matrices to subspaces, see \S 4. Moreover, a recent paper \cite{steini2} establishes a connection with integrable systems and Hilbert transform identities. We note that the problem has a natural physical interpretation referred to as Gauss' electrostatic interpretation, based on the identity
$$ \frac{p_n'(x)}{p_n(x)} = \sum_{k=1}^{n}{ \frac{1}{x-x_k}}.$$
In particular, roots of  the polynomial $p_n$ may be interpreted as point charges confined to a line, exerting a mutually repulsive force. The critical points, which coincide with the roots of $p_n'$, are equilibrium points at which these forces cancel. 
There are now three existing approaches to this open problem.
\begin{enumerate}
\item \textbf{The case $t=0$.} For random roots in the complex plane, it was conjectured by Pemantle \& Rivin \cite{pem} that the roots of $p_n'$ are distributed according to the same measure as the roots of $p_n$. By induction, this holds for any fixed derivative $p_n^{(k)}$. This was proven by Kabluchko \cite{kab}, we also refer to  \cite{cheung, red, sub}. O'Rourke \& Williams \cite{or, or3} and Kabluchko \& Seidel \cite{kab2} have provided a very fine analysis. In the simple special case where all roots are real-valued, the interlacing phenomenon immediately implies these results.
\item \textbf{The case $0 < t < 1$.} Much less is known in this case.  In \cite{steini} it was suggested that if the roots are distributed according to a smooth density $u(0,x)$ such that $\left\{x \in \mathbb{R}: u(0,x) > 0 \right\}$ is a finite interval, then the density $u(t,x)$ of roots of the $(t\cdot n)-$th derivative may be governed by the nonlinear and nonlocal evolution equation
$$ \frac{\partial u}{\partial t} + \frac{1}{\pi} \frac{\partial}{\partial x}\arctan{ \left( \frac{Hu}{ u}\right)}  = 0 \qquad \mbox{on}~\left\{x: u(x) > 0\right\},$$
where 
$$ Hf(x) =  \mbox{p.v.}\frac{1}{\pi} \int_{\mathbb{R}}{\frac{f(y)}{x-y} dy} \qquad \mbox{is the Hilbert transform.}$$
The  above partial differential equation correctly predicts the behavior of the roots of derivatives for Hermite polynomials (where $u(t,x)$ follows a semicircle distribution) and Laguerre polynomials (where $u(t,x)$ is given by a Marchenko-Pastur distribution). The derivation was, however, based on heuristic arguments, a 
rigorous proof is still outstanding. It was recently shown \cite{steini2} that there are an infinite number of conservation laws that are satisfied by both explicit closed form solutions, indicating that the PDE might have an abundance of interesting structure (see also \cite{granero, or2}).
\item \textbf{The case $t=1$.} This is the case discussed in this paper.
\end{enumerate}
\subsection{Related results.}
There are a large number of results pertaining to the roots of $p_n'$ and related properties (see \cite{bruj, bruj2, branko, dimitrov, gauss,  han,  kab, kab2, lucas, malamud, or, or3, or4, pem, rahman, red, riesz, sub, sz0, sz, totik, ull, van}) as well as to the roots of higher derivatives (see \cite{byun, cheung, ravi}). Additionally, the analogous `infinite degree' setting, in which polynomials are replaced by analytic functions, has also been considered - see Polya \cite{pol}, Farmer \& Rhoades \cite{farmer} and Pemantle \& Subramanian \cite{pem2}. 
Another natural extension is to consider the dynamics of the roots of successive derivatives for complex-valued polynomials $p_n:\mathbb{C} \rightarrow \mathbb{C}.$ If the roots are given by a smooth probability distribution $u(0,z):\mathbb{C} \rightarrow \mathbb{R}_{\geq 0}$ and the limiting measure $u(0,z)$ is radial, O'Rourke \& Steinerberger \cite{or2} suggest that the following nonlocal transport equation
$$ \frac{\partial \psi}{\partial t} = \frac{\partial}{\partial x} \left( \left( \frac{1}{x} \int_{0}^{x} \psi(s) ds \right)^{-1} \psi(x) \right)$$
describes the evolution of the radial profile. This evolution equation correctly predicts the evolution of random Taylor polynomials; it would be interesting to have 
a better understanding of the behavior of this equation.

\section{Main Result}
We are motivated by trying to understand the dynamical evolution of the roots of the $(t \cdot n)-$th derivative of $p_n$. For this purpose, we have developed a fast and efficient numerical algorithm that allows us to track the evolution of the distributions. This numerical algorithm is outlined in \S 3, results obtained via the algorithm can be found throughout the paper.
We find that the underlying process does indeed seem to have a well-defined evolution. Moreover, the process seems to be regularizing (which if true, would answer a question of Polya \cite{pol}, we refer specifically to Farmer \& Rhoades \cite{farmer} and also the discussion in \cite{steini2}). Moreover, the solutions of the process close to $t=1$ seem to be characterized by a semicircle density (after possibly removing the symmetry under dilations). It is easy to see that the two known closed-form solutions, the semicircle solution and the Marchenko-Pastur solution \cite{steini}, do indeed have an asymptotic semicircle profile at $t=1$.

   \begin{figure}[H]
    \centering
    \begin{subfigure}[b]{0.49\textwidth}
        \includegraphics[width=0.97\textwidth]{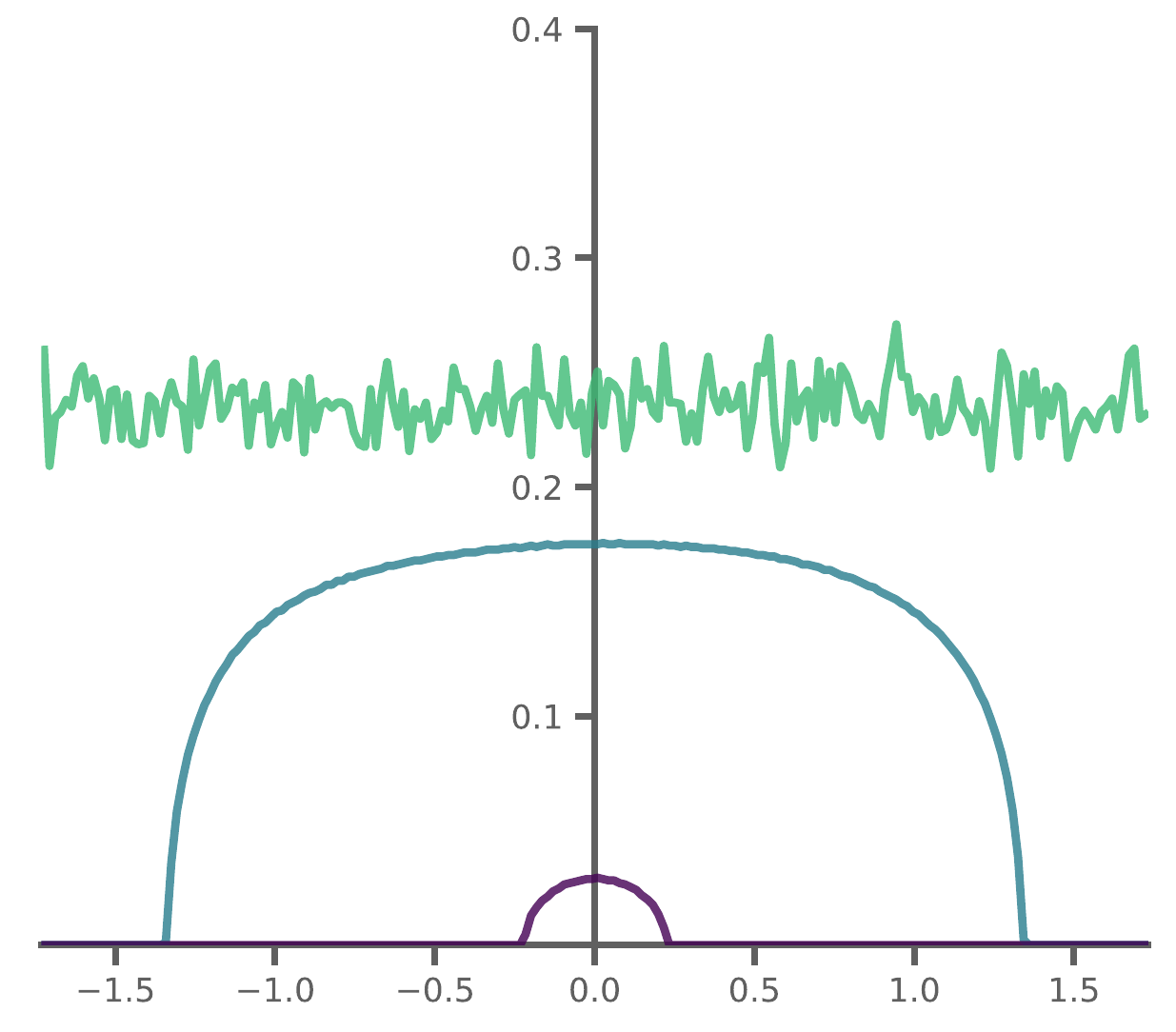}
        \label{fig:flat}
    \end{subfigure}
    ~ 
    \begin{subfigure}[b]{0.49\textwidth}
        \includegraphics[width=0.97\textwidth]{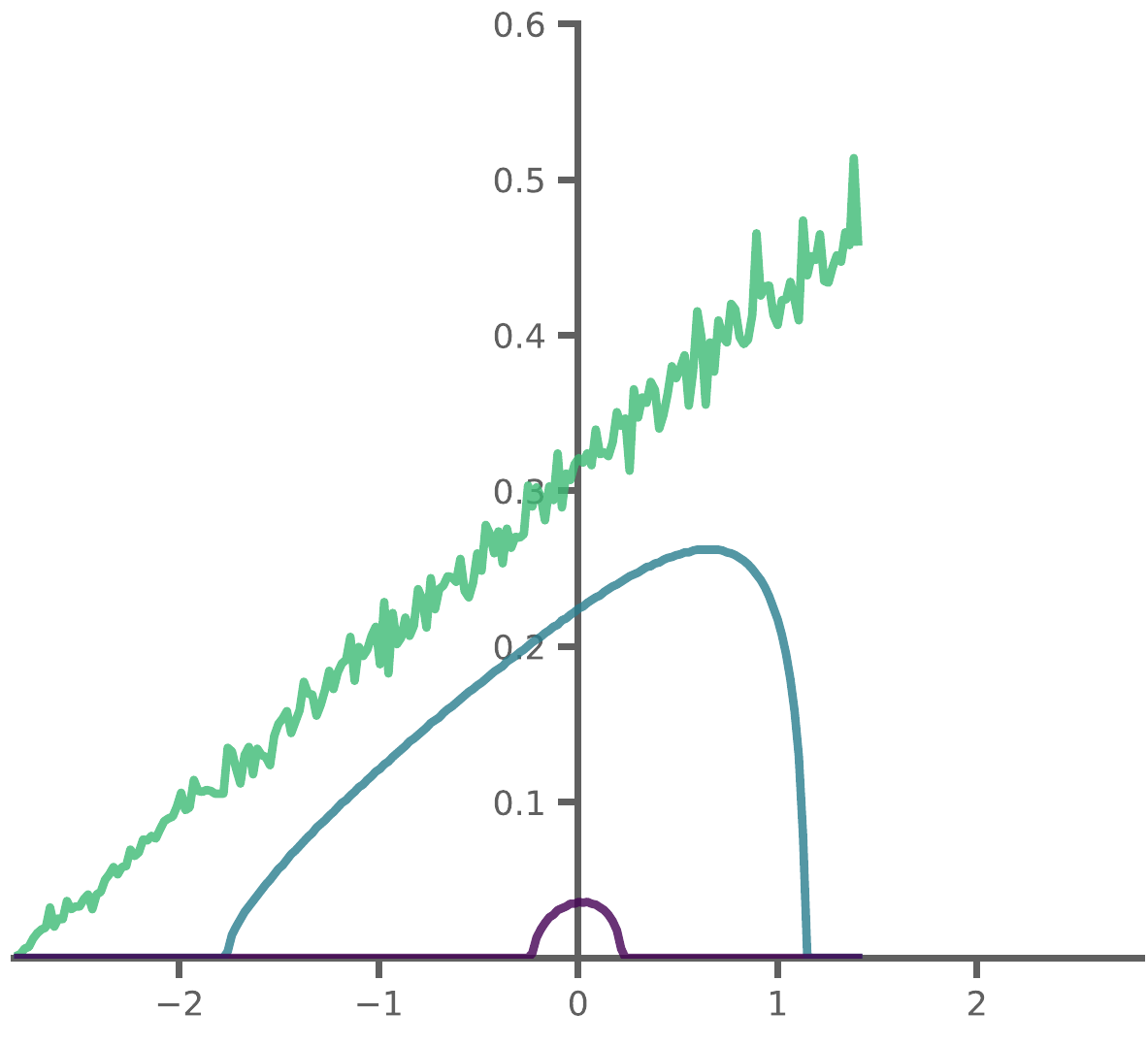}
    \end{subfigure}\\\vspace{5 ex}
    ~ 
    \begin{subfigure}[b]{0.49\textwidth}
        \includegraphics[width=0.97\textwidth]{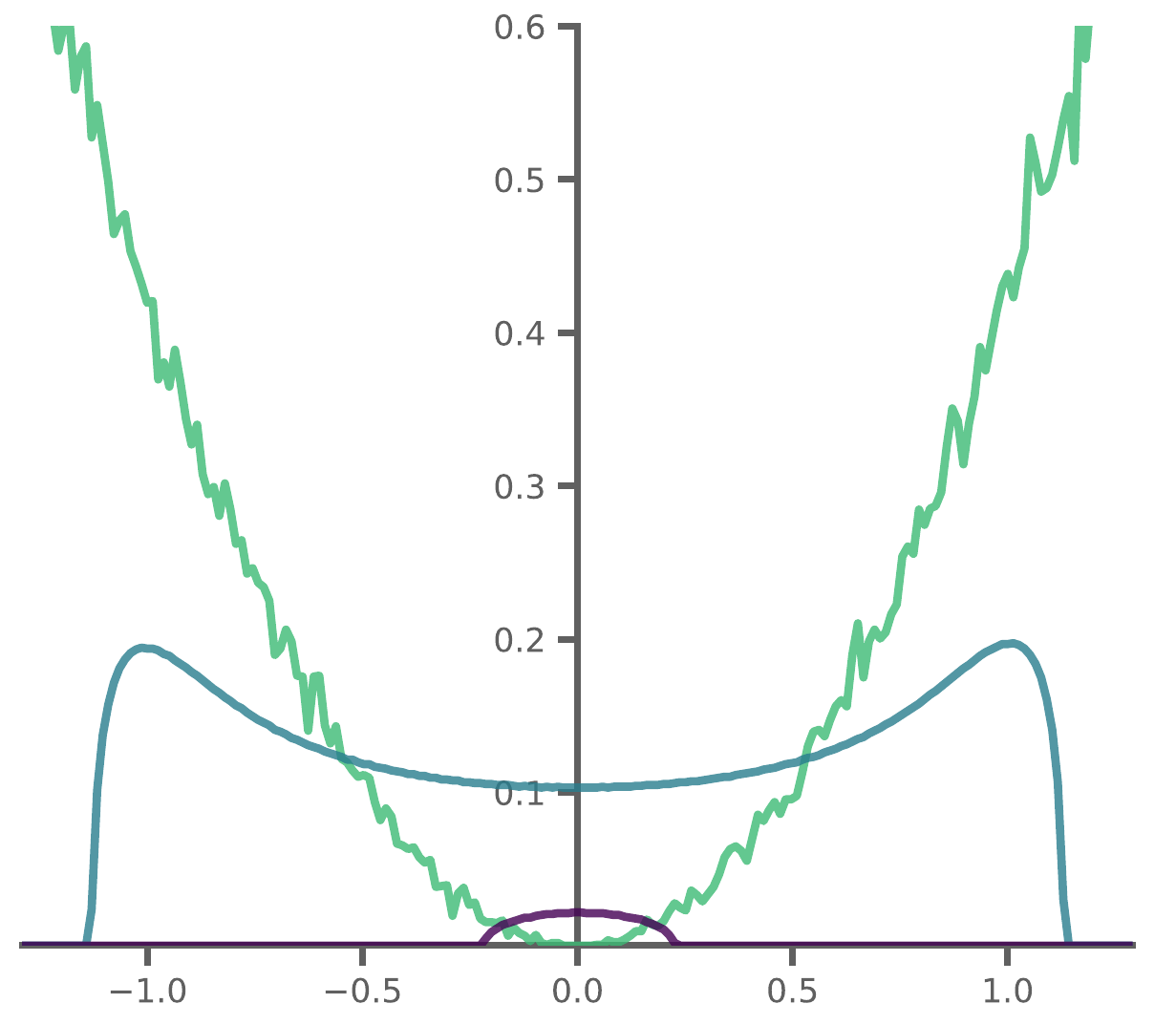}
        \label{fig:parab}
    \end{subfigure}
            \begin{subfigure}[b]{0.49\textwidth}
        \includegraphics[width=0.97\textwidth]{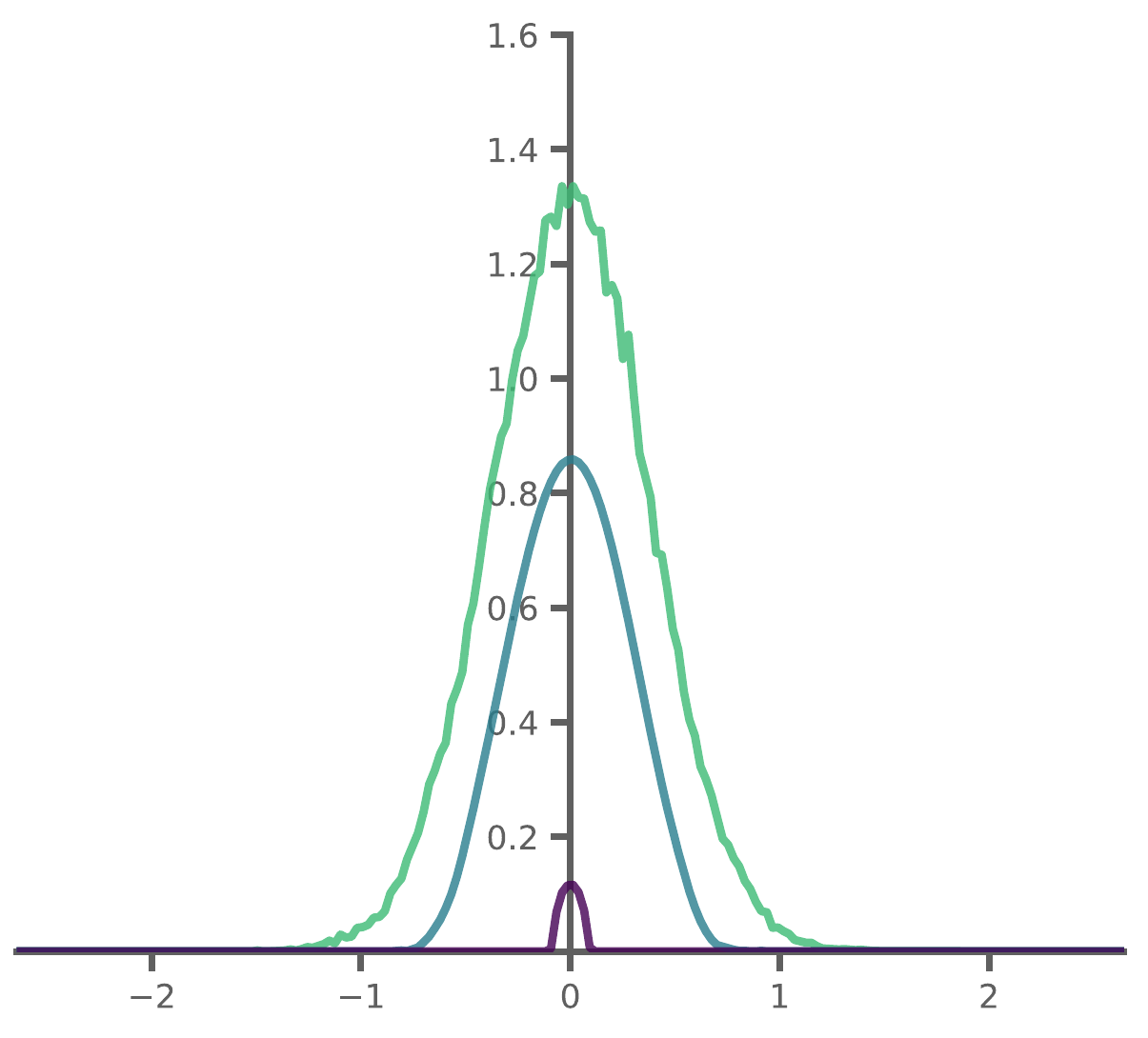}
        \label{fig:gauss}
    \end{subfigure}
    \caption{Histograms showing the evolution of roots. The green curves show the initial distribution of $80000$ roots, the blue curves show the roots after $40000$ differentiations, and the purple curves show the roots after $79000$ differentiations. All these distributions look like a semicircle at the end (up to dilation). }\label{fig:evolutions}
\end{figure}

Before stating the phenomenon, we fix some notation. We assume that we are given a probability measure on the real line such that  $\mathbb{E} |X|^k < \infty$ for all $k \in \mathbb{N}$.
Since the process under consideration is invariant under translations and dilations we can normalize it to $\mathbb{E}X = 0$ and $\mathbb{E} X^2 =1$ without loss of generality. We recall that the probabilists' Hermite polynomials are defined via the following formula
$$ He_n = (-1)^n e^{x^2/2} \frac{d^n}{dx^n} e^{-x^2/2}.$$
They deviate from the classical Hermite polynomials by a scaling factor. The first few probabilists' Hermite polynomials are given by $He_0(x) = 1$, 
$He_1(x) = x$ and
\begin{align*}
He_2(x) &= x^2-1\\
He_3(x) &= x^3 - 3x\\
He_4(x) &= x^4 - 6x^2 + 3
\end{align*}
 Their roots are known to have a Wigner semicircle distribution (up to dilation) as $n \rightarrow \infty$, see, for example, Kornik \& Michaletzky \cite{kornik}. In the following we assume that $p_n$ is a random $n$-th degree polynomial whose roots are i.i.d. random variables with a distribution satisfying  $\mathbb{E}X = 0,$ $\mathbb{E} X^2 =1$ and $\mathbb{E} |X|^k < \infty$ for all $k \in \mathbb{N}$.

\begin{thm}
Fix  $\ell \in \mathbb{N}$ and a compact interval $I \subset \mathbb{R}$. The $(n-\ell)-$th derivative of $p_n$ satisfies, for all $x \in I$, as $n \rightarrow \infty$,
 $${n^{\ell/2}} \frac{\ell!}{n!} \cdot p_n^{(n-\ell)}\left( \frac{x}{\sqrt{n}}\right) = (1+ o(1)) \cdot He_{\ell}(x + \gamma_n),$$
where $He_{\ell}$ is the $\ell-$th probabilists' Hermite polynomial and $\gamma_n$ is random variable converging to the standard $\mathcal{N}(0,1)$ Gaussian as $n \rightarrow \infty$. 
\end{thm}

\begin{center}
\begin{figure}[h!]
        \includegraphics[width=0.6\textwidth]{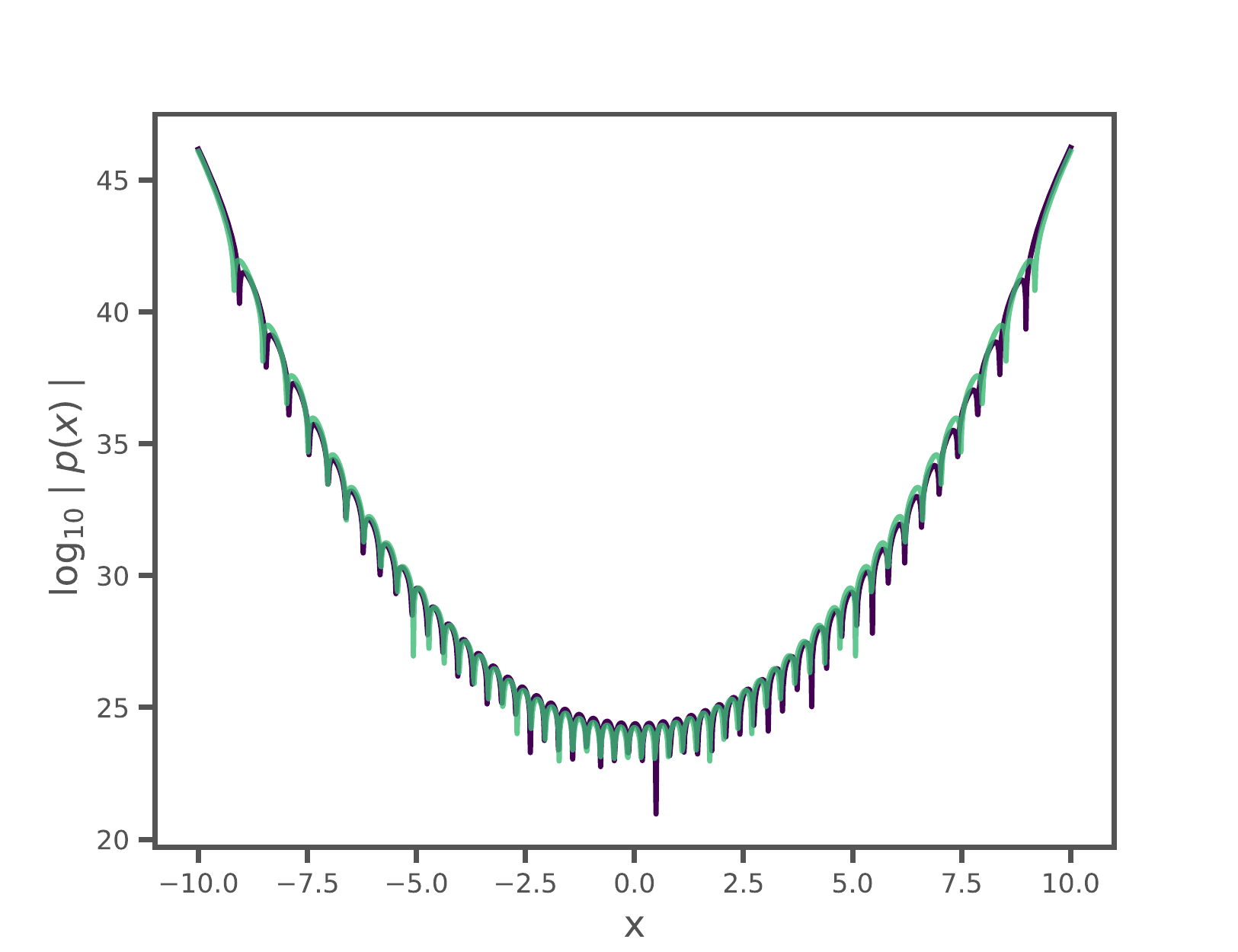}
        \caption{The Theorem illustrated: after computing 950 derivatives of a random polynomial of degree 1000, the resulting polynomial, after a suitable shift and scaling is closely matched by $He_{50}$ (both shown on logarithmic scale).}
\end{figure}
\end{center}

 $n^{\ell/2} \ell!/n!$ is merely a prefactor and has no influence on the roots. The shift by a Gaussian is to be expected: it is an elementary fact (see \cite{rahman}) that the means of the roots of $p_n$ and $p_n'$ coincide. This is also the second conservation law in \cite{steini2}. The mean of $n$ random variables having mean 0 is an approximately Gaussian random variable $\sim n^{-1/2} \cdot \mathcal{N}(0,1)$. After differentiating $n-\ell$ times, the remaining $\ell$ roots will be shown in \S 4 to be at scale $\sim n^{-1/2}$ and this explains the shift. In particular, the case $\ell = 1$ is equivalent to the central limit theorem.\\
 
 \textbf{Open Problems.}
 This result motivates a number of interesting problems. Maybe the most natural question is whether it is possible to let $\ell$ grow with $n$.  It is conceivable that our proof could be adapted to show that it is possible to have, say, $\ell \sim \log\log{n}$ but it is not clear to us whether $\ell$ can maybe even be as large as $\ell \sim n/\log \log n$. It is clear from the Marchenko-Pastur solution (see \cite{steini}) that $\ell$ cannot grow linearly in $n$. It would also be interesting to understand what happens when $\ell \sim \varepsilon n$: though the distribution will not be exactly a semicircle, our numerical examples show that in most instances it is quite close. How close is it and how does this depend on the initial probability distribution?

  \section{The Algorithm}
  This section discusses the algorithm used to investigate the dynamical evolution and produce the figures in this paper. We also discuss some intuitiongained from these figures.
  
  \subsection{Typical Evolution.} 
  The evolution of roots under repeated differentiation is shown in Figure \ref{fig:evolutions} for a representative set of initial distributions. Empirically, we observe that this process is indeed smoothing: inhomogeneities in the density even out and mass moves from regions with high density to those with lower density. This is in line with 
 known conjectures \cite{farmer, pol, steini2}. Moreover, the observed behaviour supports an interesting hypothesis about the PDE obtained in \cite{steini}: does the nonlocal evolution equation
 $$ \frac{\partial u}{\partial t} + \frac{1}{\pi} \frac{\partial}{\partial x}\arctan{ \left( \frac{Hu}{ u}\right)}  = 0 \qquad \mbox{on}~\left\{x: u(x) > 0\right\}$$
increase the regularity of its solution? Transport equations usually have no reason to increase regularity, however, this equation is also driven by a nonlocal term that may have a positive impact on the regularity.

  \subsection{The Gap Filling.} One particularly interesting open problem is to understand the gap filling mechanism: if the initial distribution of roots has a `gap' - an internal interval on which it vanishes, then successive differentiation will introduce roots into the initially-empty interval. It is not clear how the density of these roots depends on the initial densities to the left and right of the gap. 
  
    \begin{figure}[h!]
\centering 
\begin{tikzpicture}[scale=1]
\draw [thick] (0,0) -- (6,0); 
\draw [thick] (1,0) to[out=60, in=110] (3,0);
\draw [thick] (4,0) to[out=80, in=110] (5,0);
\node at (3.5, 0.3) {$?$};
\end{tikzpicture}
\caption{A density not covered by the analysis in \cite{steini}.}
\end{figure}
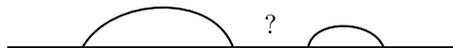
  
\begin{figure}[H]
    \centering
    \begin{subfigure}[b]{0.49\textwidth}
        \includegraphics[width=0.99\textwidth]{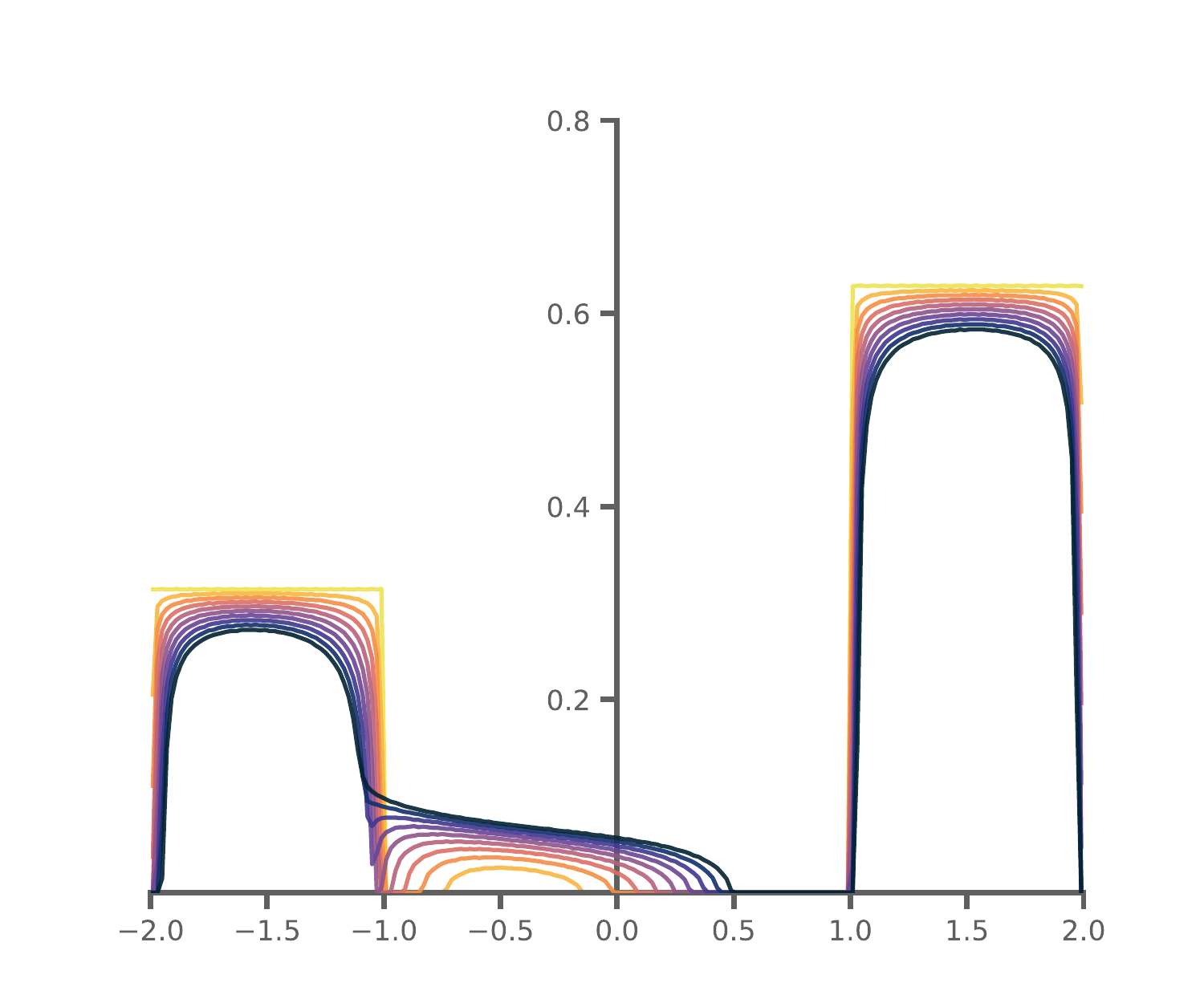}
        \label{fig:first_gap}
    \end{subfigure}
      \begin{subfigure}[b]{0.49\textwidth}
        \includegraphics[width=0.99\textwidth]{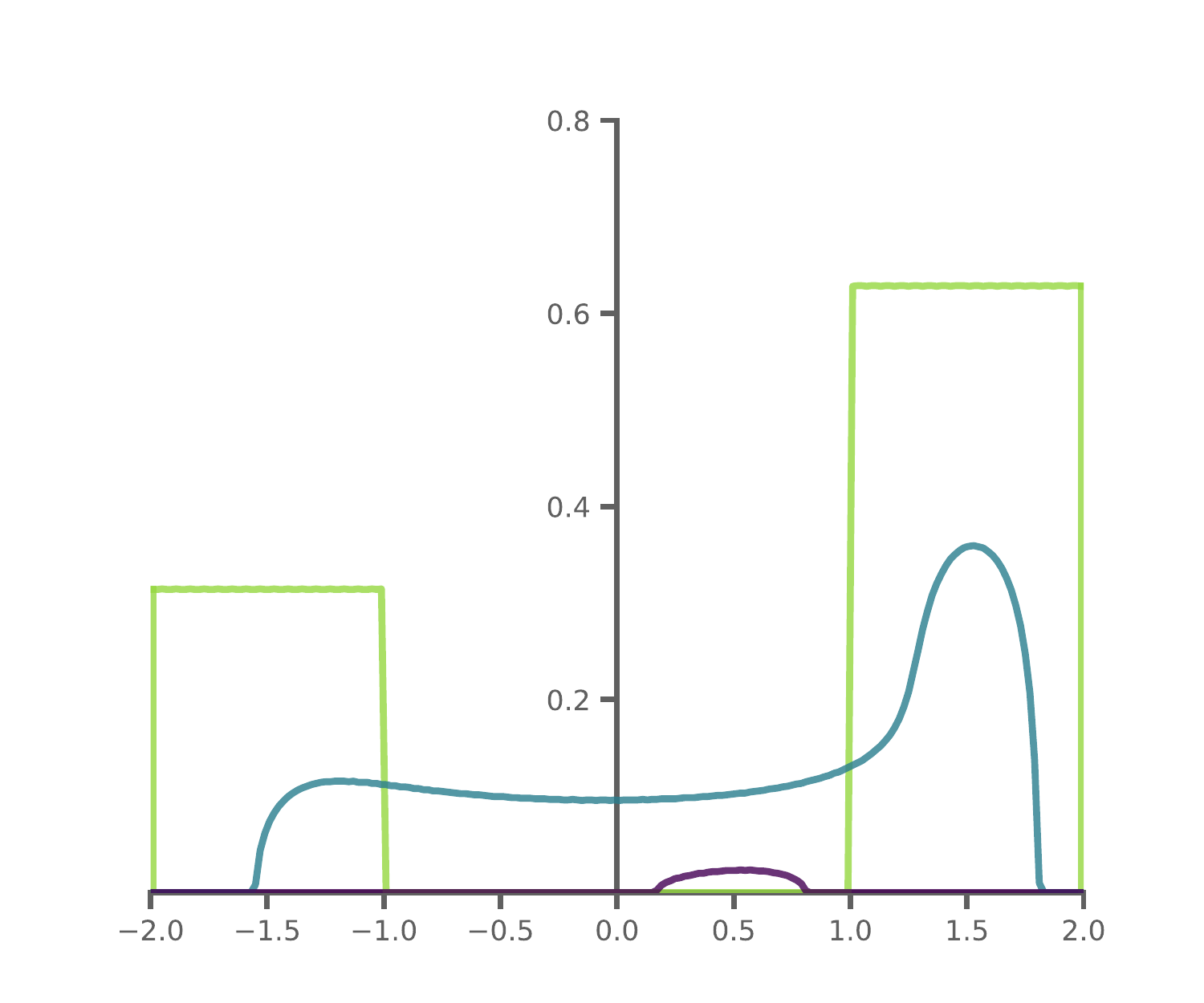}
        \label{fig:first_gap}
    \end{subfigure}
    \caption{Histograms showing the evolution of roots for an initial distribution with a gap. left: from light to dark the roots after 0, 1000, 2000, 3000, 4000, 5000, 6000, 7000, 8000 and 9000 derivatives. right: roots for the same initial distribution after 0, 40000, and 79000 differentiations. The initial number of roots was 80000.}\label{fig:evolutions_gap}
\end{figure}

  In particular, the derivation of the PDE in \cite{steini, steini2} requires that the density be supported on a connected set on which it does not vanish. Currently there are no predictions on how the roots of the $(t\cdot n)-$th derivative of $p_n$ evolves when the roots of $p_n$ are, say, distributed according to the probability measure
  $$ \mu = \frac{1}{3} \chi_{[-2,-1]} + \frac{2}{3} \chi_{[1,2]}.$$
  The algorithm implemented in this paper may serve as a useful tool for obtaining intuition about the behavior in this context. The evolution for the above density is shown in Figure \ref{fig:evolutions_gap}.
  
 We observe something quite interesting: though the roots slowly fill the interval, in contrast to, say, parabolic equations, this process is not instantaneous. A careful inspection of
$$ \frac{p_n'(x)}{p_n(x)} = \sum_{k=1}^{n}{\frac{1}{x-x_k}}$$
shows that, for this particular example, roots on the left end in the right bump do indeed initially move towards the right and that no roots are being created in $[0.9, 1]$ for quite some time. Numerical experiments indicate a wealth of structure. A particularly interesting example is given by the initial distribution
$$ \mu = \frac{n}{2} \delta_{-1} + \frac{n}{2} \delta_{1}.$$
We observe empirically that upon differentiation, roots are being created close to the origin and that the initial profile seems to be that of a semicircle (at least in the infinitesimal sense when the number of derivatives is a fixed integer $\ell$ and $n \rightarrow \infty$). We can rigorously establish this phenomenon.
  
  \begin{proposition} For $\ell \in \mathbb{N}$ fixed as $n \rightarrow \infty$, the $\ell-$th derivative of the polynomial $p(x) = (x^2 -1)^n$ evaluated at scale $y=x/\sqrt{n}$ converge to the $\ell-$th Hermite polynomial. In particular, the distribution of roots forms a semicircle.
  \end{proposition}
  
  \subsection{The Idea behind the Algorithm.} Given a polynomial $p_n$ of degree $n$ with $n$ distinct roots $x_1,\dots,x_n,$ a real number $r$ is a root of its derivative $p_n'$ if and only if it satisfies the equation
  \begin{equation*}\label{eqn:root_sum}
  \sum_{i=1}^{n} \frac{1}{r-x_i} = 0.
  \end{equation*}
  Moreover, by the interlacing property the above equation has exactly one zero on each interval $(x_i,x_{i+1}),$ $i=1,\dots,n-1.$ Our algorithm obtains the roots of $p_n'$ sequentially by applying Newton's method in each of these intervals.
  It is easily seen that a naive implementation of this approach would require $O(n^2)$ operations to obtain all $n-1$ roots of $p_n'.$ When the number of roots is large and many derivatives are required this cost may become prohibitively expensive. On the other hand, the computations of sums of this form arise frequently in physics and are amenable to a variety of fast algorithms. We use a minor modification of the algorithm proposed in \cite{gimb}, which is well-adapted to this context. \\
  
  Given a point $r$ at which the sum is to be evaluated, the method splits the sum into {\it local} contributions (those associated with roots that are within some pre-defined distance of $r$) and {\it farfield} contributions (the portion of the sum due to the remainder of the roots). The computation of the local contribution is done by directly evaluating the corresponding terms in the sum, and requires $\mathcal{O}(1)$ operations per choice of $r.$ Using the method in \cite{gimb}, for any prescribed accuracy $\varepsilon>0$ the farfield contribution of $x_1,\dots,x_n$ evaluated at $k$ points can be computed in $O(k+n\log (k\, \varepsilon^{-1}))$ operations. 
  Hence, evaluating
  $$\sum_{i=1}^n \frac{1}{r-x_i}$$
  at $\mathcal{O}(n)$ points can be performed in $\mathcal{O}(n \log{(n \varepsilon^{-1})})$ operations. Unlike in \cite{gimb}, rather than evaluating the farfield contributions at the original points $x_1,\dots,x_n,$ we evaluate them at a set of interpolation nodes which can then be used to calculate the farfield contribution at any point in $\mathcal{O}(1)$ operations after a precomputation step requiring $\mathcal{O}(n \log{(n \varepsilon^{-1})})$ operations. Alternate approaches based on slightly different fast algorithms for farfield computations can be found in \cite{dutt,gu,greengard} and the references therein, for example. We refer to \cite{pan,schleicher} and the references therein for efficient algorithms for the related problem of efficient polynomial root-finding.\\

  \textbf{Runtime and Accuracy.} The algorithm was implemented in gfortran and run on a 2019 MacBook Pro with 16 GB of memory and a 2.5 GHz processor. As an illustration of its speed and accuracy it was run on the 10000${}^{\rm th}$ Hermite polynomial $He_{10000}(x).$ Analytically, after 5000 differentiations the roots correspond to those of the 5000${}^{\rm th}$ Hermite polynomial. The maximum error in the roots of $He_{5000}$ computed by successive differentiation was $3.5527 \times 10^{-13}$ and it took 54 seconds to compute the 37,497,500 roots of $He_{9999},\dots,He_{5000}.$  

  \section{Proofs}
\subsection{Preliminaries}  We use the notation $e_k(x_1, \dots, x_n)$, which we abbreviate as $e_k$, to denote the $k-$th elementary symmetric polynomial on $n$ variables, i.e.
\begin{align*}
e_0(x_1, \dots, x_n) &= 1\\
e_1(x_1, \dots, x_n) &= x_1 + \dots + x_{n}\\
e_2(x_1, \dots, x_n) &= \sum_{i < j}{x_i x_j}\\
e_3(x_1, \dots, x_n) &= \sum_{i < j < k}{x_i x_j x_k}
\end{align*}
and so on for $k \leq n$. For $k>n$, we set $e_k(x_1, \dots, x_n) = 0$. We define the $k-$th power sum as 
$$ p_k(x_1, \dots, x_n) = x_1^k + x_2^k + \dots + x_n^k$$
and will abbreviate it as $p_k$.
The elementary symmetric polynomials arise naturally as the coefficients in the monic polynomial determined by the roots $x_1, \dots, x_n$
$$ \prod_{i=1}^{n}{(x-x_i)} = \sum_{k=0}^{n} (-1)^k e_k(x_1, \dots, x_n) x^{n-k}.$$
We differentiate this polynomial $(n-\ell)-$times to arrive at
$$  p_n^{(n-\ell)}(x) = \sum_{k=0}^{\ell} (-1)^k e_k x^{\ell-k} \frac{(n-k)!}{(l-k)!}.$$
We normalize this polynomial to be monic and obtain
$$ \frac{\ell!}{n!} \cdot p_n^{(n-\ell)}(x) = \sum_{k=0}^{\ell} (-1)^k e_k x^{\ell-k} \frac{(n-k)!}{(l-k)!} \frac{\ell!}{n!}.$$
 Introducing new variables
$$ f_k = e_k  \frac{\ell!}{(\ell-k)!} \frac{(n-k)!}{n!},$$
we can write this as
$$ \frac{\ell!}{n!}  \cdot p_n^{(n-\ell)}(x) = \prod_{k=1}^{\ell}{(x - r_i)} = \sum_{k=0}^{\ell} (-1)^k f_k x^{\ell-k},$$
where $r_1, \dots, r_{\ell}$ denote the roots of the $(n-\ell)-$th derivative of $p_n$ and $f_k$ is
the $k-$th elementary symmetric polynomial applied to these roots $r_1, \dots, r_{\ell}$.

\subsection{Mean and Variance.} In this section, we use the Newton identities to obtain estimates on the mean
and variance of the $\ell$ roots of the $(n-\ell)-$th derivative of $p_n$. This computation is not required to
prove the Theorem but is useful in building some intuition and in explaining the parameter choices in
the proof.\\

\textit{Mean.} Since $x_1, \dots, x_n$ are i.i.d. random variables with mean value 0 and variance 1, we have a fairly good
understanding of the mean
$$ \overline{x} = \frac{x_1 + \dots + x_n}{n}$$
and expect it to be roughly distributed like $\sim n^{-1/2} \cdot \mathcal{N}(0,1)$. We will now show that this mean is invariant under differentiation. 
By writing
$$ p_n(x) = \prod_{i=1}^{n}{ (x-x_i)} =  x^n + \sum_{k=0}^{n-1}{ a_k x^k},$$
we see that
$$ \sum_{i=1}^{n}{x_i} = - a_{n-1}.$$
Let us now consider the $(n-\ell)-$th derivative of the polynomial $p_n(x)$
$$ p_n^{(n-\ell)}(x) = \frac{n!}{\ell!} x^{\ell} + a_{n-1} \frac{(n-1)!}{(\ell -1)!} x^{\ell -1} + \dots =\frac{n!}{\ell!} \cdot \prod_{i=1}^{\ell}{ (x - r_i)}.$$
Normalizing the polynomial to be monic (which has no effect on the roots) and computing the coefficients again shows that
$$ \sum_{i=1}^{\ell}{r_i} = - a_{n-1} \frac{\ell}{n}$$
and thus
$$ \overline{r} =  \frac{1}{\ell}  \sum_{i=1}^{\ell}{r_i} = -\frac{a_{n-1}}{n} = \overline{x}.$$
In particular, since $\overline{x}$ behaves like $\sim n^{-1/2}\cdot \mathcal{N}(0,1)$, the same will be true for the mean of the roots of the $(n-\ell)-$th derivative. This accounts for our rescaling by $\sim n^{-1/2}$ and the arising shift by a random Gaussian variable.\\

\textit{Variance.} For the variance, we use the following identity which can be found in the book of Rahman \& Schmeisser \cite[Lemma 6.1.5]{rahman}: using $x_1, \dots, x_n$ to denote the roots of the polynomial $p_n$ and $r_1, \dots, r_{\ell}$ to denote the roots of the polynomial $p_n^{(n-\ell)}$, we have the identity
$$ \frac{1}{n^2(n-1)} \sum_{1 \leq i < j \leq n}{ (x_i - x_j)^2} = \frac{1}{\ell^2(\ell-1)} \sum_{1 \leq i < j \leq \ell}{ (r_i - r_j)^2}.$$
We refer to \cite{steini2} for a more in-depth discussion of such conservation laws. We can rewrite this slightly more informally as
$$ \mbox{average of} ~(r_i - r_j)^2 = \frac{\ell}{n}\cdot \mbox{average of}~ (x_i - x_j)^2.$$
In particular, since the variance decreases, the roots of the derivative move closer to each other. By assumption
$$ \mbox{average of}~ (x_i - x_j)^2 = \mathbb{E} (x_i - x_j)^2 = 2,$$
and so the roots of the $(n-\ell)-th$ derivative satisfy the identity
$$ \mbox{average of} ~(r_i - r_j)^2 = \frac{2 \ell}{n},$$
from which it follows that a typical root $r_i$ is roughly at distance $\sim \ell^{1/2} n^{-1/2}$ from the origin.

\subsection{The Probabilistic Lemma.} 
\begin{figure}[h!]
    \centering
        \includegraphics[width=0.7\textwidth]{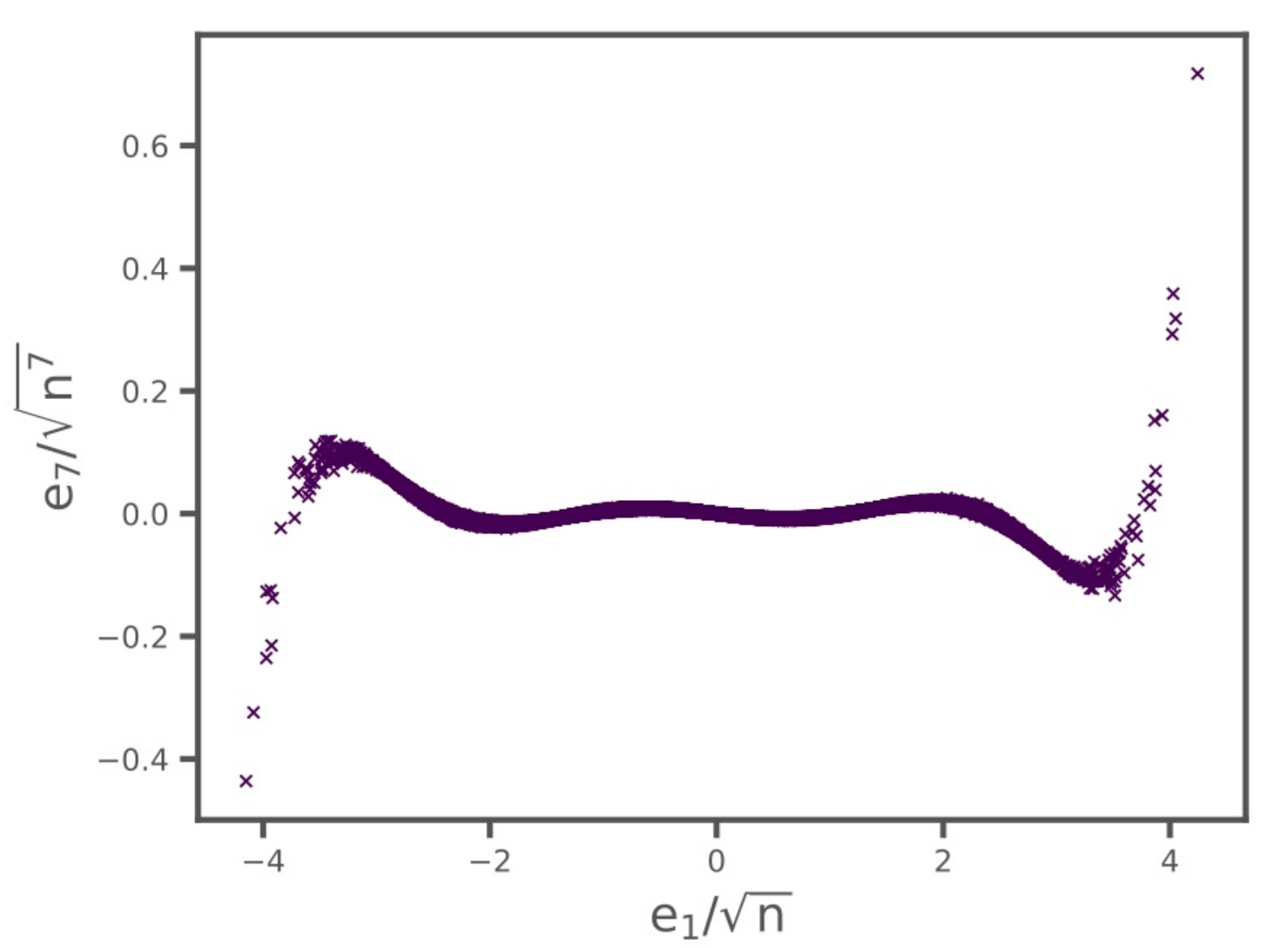}
    \caption{Scaled values of $e_1$ and $e_7$ calculated using $n=1000$ points chosen uniformly on the interval $[-\sqrt{3},\sqrt{3}].$}\label{fig:esums}
\end{figure}

This section contains our main probabilistic ingredient: a concentration result for
elementary symmetric polynomials. In short, it says that if $x_1, \dots, x_n$ are i.i.d. random variables coming from a distribution
on the real line with $\mathbb{E}X = 0$, $\mathbb{E}~X^2 = 1$ and finite moments $\mathbb{E} |X|^k < \infty$, then the elementary symmetric polynomial $e_k$ is, to leading order, given
by a polynomial of $e_1$ and $n$ (as long as $k$ is small compared to $n$) -- see Figure \ref{fig:esums}. One way of phrasing this result is that
$$ e_1(x_1, \dots, x_n) \qquad \mbox{and} \qquad e_k(x_1, \dots, x_n)$$
are tightly correlated in the sense of Chatterjee \cite{cha} for any fixed $k$ as $n \rightarrow \infty$.
  
We illustrated the first few cases of the main result to illustrate the idea. It will also serve as the base case for the induction step in the proof of the Lemma.
The proof uses the Newton identities
$$ m \cdot e_m(x_1, \dots, x_n) = \sum_{i=1}^{m} (-1)^{i-1} e_{m-i}(x_1, \dots, x_n) p_i(x_1, \dots, x_n).$$
In particular, for $m=2$ we have
$$ e_2 = \frac{e_1^2}{2} - \sum_{i=1}^{n}{x_i^2}.$$
However, since the $x_i$ are i.i.d. random variables and $\mathbb{E}X = 0$ as well as $\mathbb{E}X^2 = 1$, we obtain
$$ e_2 = \frac{e_1^2}{2} - n + \mbox{error},$$
where the error behaves like $\mbox{error} \sim \sqrt{n} \cdot \mathcal{N}(0,1)$.
This shows that $e_2$ can be expressed to leading order as a polynomial of $e_1$ and $n$.
The Newton identities coupled with our existing result show that
\begin{align*}
e_3 &= \frac{e_2}{3} e_1 - \frac{e_1}{3} \sum_{i=1}^{n}{x_i^2}  + \frac{1}{3} \sum_{i=1}^{n}{x_i^3}\\
&= \frac{1}{3} \left( \frac{e_1^2}{2} - n + \mbox{error}\right)e_1 - \frac{e_1}{3} \left(n + \mbox{error}\right) + \frac{n}{3}\mathbb{E} X^3 + \mbox{error}\\
&= \frac{e_1^3 - 5n e_1}{6} + \mbox{error}_2
\end{align*}
where each $\mbox{error}$ term is comparable to $\sqrt{n} \cdot \mathcal{N}(0,1)$ making $\mbox{error}_2$ a lower order term comparable to $\sim n \cdot \mathcal{N}(0,1)$. Again, we observe that $e_3$ is a polynomial in terms of $e_1$ and $n$ up to a lower order term. We now state the general result.

\begin{lem} Let $m \in \mathbb{N}$ and let $x_1, \dots, x_n$ be i.i.d. random variables sampled from a distribution on $\mathbb{R}$ with $\mathbb{E}X = 0$, $\mathbb{E} X^2 = 1$ and $\mathbb{E} |X|^m < \infty$. Then, as $n \rightarrow \infty$,
 $$\mathbb{E} ~\left|  e_m -  \sum_{k=0}^{\left\lfloor m/2 \right\rfloor} (-1)^k \frac{1}{k! (m-2k)! 2^k} \cdot e_1^{m-2k} n^{k} \right| \lesssim_{X} n^{\frac{m-1}{2}}$$ 
\end{lem}
We note that we expect $e_m$ to be at scale $\sim n^{m/2}$ which shows that this sum captures the main contribution up to an order at a smaller scale. The implicit constant could be estimated in terms of the growth of the first few moments. We recall that this sum is related to the probabilists Hermite polynomials
$$ H e_m(x) =  \sum_{k=0}^{\left\lfloor m/2 \right\rfloor} \frac{m!}{k! (m-2k)! 2^k }  \cdot x^{m-2k}.$$
 This relates the Lemma to results of Rubin \& Vitale \cite{rubin} and Mori \& Szekely \cite{mori}. Our argument seems different and is fairly elementary.

\begin{proof}[Proof of the Lemma]
The proof proceeds via induction on $m$. We have already seen the cases $m=1, m=2$ and $m=3$. As for the general case, we note that since $e_m$ is
a sum over products of different independent random variables, we have, using the Cauchy-Schwarz inequality and the fact that $\mathbb{E}X = 0$,
\begin{align*}
\left(\mathbb{E}~ |e_m|\right)^2 \leq \mathbb{E} ~e_m^2 = \binom{n}{m} \sim c_m n^m
\end{align*}
and we thus expect $e_m$ to be typically at scale $\sim n^{m/2}$.
Let us now assume the statement has been verified up to $m-1$. The Newton identities imply
$$ m \cdot e_m = \sum_{i=1}^{m} (-1)^{i-1} e_{m-i}  \cdot \sum_{j=1}^{n}{x_j^i}.$$
We observe that all the power sums are actually relatively small with high likelihood: since the moments are finite, we expect
$$  \mathbb{E} \left| \sum_{j=1}^{n}{x_j^i} \right| \lesssim  n \cdot \mathbb{E} |X|^i \lesssim_{X} n.$$
This leads us to suspect that only the first two terms in the Newton expansion are actually relevant. Indeed, we have
$$ e_m - \frac{e_{m-1} e_1}{m} - \frac{e_{m-2}}{m} \sum_{j=1}^{n}{x_j^2} = \sum_{i=3}^{m} (-1)^{i-1} e_{m-i}  \cdot \sum_{j=1}^{n}{x_j^i}.$$
Taking the absolute value and expectation on both sides, we see that 
\begin{align*}
\mathbb{E} ~\left|e_m - \frac{e_{m-1} e_1}{m} - \frac{e_{m-2}}{m} \sum_{j=1}^{n}{x_j^2}  \right| &\leq \sum_{i=3}^{m} \mathbb{E}  \left|e_{m-i}  \cdot \sum_{j=1}^{n}{x_j^i} \right| . \\
&\lesssim \mathbb{E} |e_{m-3}| n \lesssim n^{\frac{m-1}{2}},
\end{align*}
where the last two inequalities follow from the inequalities above in combination with standard concentration arguments. We observe that 
$$ \sum_{j=1}^{n}{x_j^2} = n \pm \sqrt{n} \cdot \mathcal{N}(0,1).$$
We use Cauchy-Schwarz to argue that
\begin{align*}
\mathbb{E} ~\left|e_{m-2} \cdot \left(\sum_{j=1}^{n}{x_j^2}  - n \right)\right| &\leq \left( \mathbb{E}~|e_{m-2}|^2 \right)^{1/2} \left( \mathbb{E}\left(\sum_{j=1}^{n}{x_j^2}  - n \right)^2\right)^{1/2} \\
&\lesssim_{X} n^{\frac{m-2}{2}} n^{1/2} = n^{\frac{m-1}{2}}.
\end{align*}
Therefore the difference between the random variables
$$ X = e_m - \frac{e_{m-1} e_1}{m} - \frac{e_{m-2}}{m} \sum_{j=1}^{n}{x_j^2}$$
and
$$ Y = e_m - \frac{e_{m-1} e_1}{m} - \frac{e_{m-2}}{m} n$$
satisfies
\begin{align*}
\left| \mathbb{E}X - \mathbb{E} Y\right| \lesssim_X n^{\frac{m-1}{2}}.
\end{align*}
We conclude the argument by determining $Y$. We note that, by induction,
\begin{align*}
 e_{m-1}e_1 - e_{m-2}n 
 &= \frac{1}{(m-1)!} \sum_{k=0}^{\left\lfloor (m-1)/2 \right\rfloor} (-1)^k \frac{(m-1)!}{k! (m-1-2k)! 2^k} \cdot e_1^{m-2k} n^{k} \\
 &-  \frac{1}{(m-2)!} \sum_{k=0}^{\left\lfloor m/2 \right\rfloor - 1} (-1)^k \frac{(m-2)!}{k! (m-2-2k)! 2^k} \cdot e_1^{m-2-2k} n^{k+1} + Z,
 \end{align*}
 where $Z$ is a random variable satisfying
  $$ \mathbb{E} |Z| \lesssim_X n^{\frac{m-2}{2}}.$$
 We observe that the leading coefficient in front of $e_1^m$ is given by
 $$ \frac{1}{(m-1)!} \frac{(m-1)!}{0! \cdot (m-1)!} e_1^{m} = \frac{e_1^m}{(m-1)!}$$
 which is consistent with the desired formula (note the multiplication with $m$ on the right-hand side). For $k \geq 1$, the coefficient $c_k$ in the expression $ c_k e_1^{m-2k}$ is 
 \begin{align*}
 c_k &= \frac{1}{(m-1)!} (-1)^k \frac{(m-1)!}{k! (m-1-2k)! 2^k} n^k \\
 &+ \frac{1}{(m-2)!} (-1)^{k} \frac{(m-2)!}{(k-1)! (m-2k)! 2^{k-1}} n^k\\
 &= (-1)^k n^k \left[ \frac{1}{k! (m-1-2k)! 2^k} + \frac{1}{(k-1)!(m-2k)! 2^{k-1}} \right] \\
 &=(-1)^k n^k \left[ \frac{m-2k}{k! (m-2k)! 2^k} + \frac{2k}{k!(m-2k)! 2^{k}} \right] =(-1)^k n^k  \frac{m}{k! (m-2k)! 2^k}.
 \end{align*}
This is the right term for the expression $m \cdot e_{m}$, by canceling the factor $m$ on both sides, we obtain the desired result by induction.
\end{proof}

\subsection{Proof of the Theorem}
\begin{proof}
We recall that, as derived above, we have
$$ S =  p_n^{(n-\ell)}(x) = \prod_{k=1}^{\ell}{(x - r_i)} =\sum_{m=0}^{\ell} (-1)^m f_m x^{\ell-m},$$
where 
$$ f_k = e_k  \frac{\ell!}{(\ell-k)!} \frac{(n-k)!}{n!} =  e_k  \frac{\ell!}{(\ell-k)!} \frac{1 + o(1)}{n^k}$$ 
Additionally, we recall that

$$  e_k = \sum_{i=0}^{\left\lfloor k/2 \right\rfloor} (-1)^i \frac{1}{i! (k-2i)! 2^i} \cdot e_1^{k-2i} n^{i} + o(n^{\frac{k}{2}}).$$
This means that the leading order term as $n \rightarrow \infty$ is given by
$$ S = \sum_{m=0}^{\ell} (-1)^m \left(  \binom{l}{m} \sum_{k=0}^{m/2} (-1)^k \frac{m!}{k! (m-2k)! 2^k} \cdot e_1^{m-2k} n^{k-m}\right) x^{\ell-m},$$
where, by a slight abuse of notation, we write $m/2$ in place of $\left\lfloor m/2 \right\rfloor$ in the limit of the summation.
We recall that $e_1 \sim n^{1/2}$ and thus
$$ e_1^{m-2k} n^{k-m} \sim n^{-m/2},$$
which in turn implies that all the summands in the inner sum are roughly comparable in size; they are all roughly at scale $n^{-m/2}$ with the lower order term being a factor $\sim n^{-1/2}$ smaller. Motivated by this reasoning, we introduce the random variable $\gamma$ by making the ansatz
$$e_1 = \gamma \sqrt{n}.$$
Recalling that 
$$ e_1 = e_1(x_1, \dots, x_n) = \sum_{i=1}^{n}{x_i},$$
The above rescaling simplifies the leading order term $S$ to
$$ S_1 = \sum_{m=0}^{\ell} (-1)^m \left(  \binom{l}{m} \sum_{k=0}^{m/2} (-1)^k \frac{m!}{k! (m-2k)! 2^k} \cdot \gamma^{m-2k} n^{-m/2}\right) x^{\ell-m}$$
The natural scale on which to evaluate this quantity, motivated by the computation of the variance above, is $x \sim n^{-1/2}$. Thus we make another substitution
$$ x = \frac{y}{\sqrt{n}}$$
which results in
$$ x^{\ell - m} = n^{m/2} n^{-\ell/2} y^{\ell - m}.$$
This shows that
$$ S_1 = \frac{1}{n^{\ell/2}}\sum_{m=0}^{\ell} (-1)^m \left(  \binom{l}{m} \sum_{k=0}^{m/2} (-1)^k \frac{m!}{k! (m-2k)! 2^k} \cdot \gamma^{m-2k} y^{\ell-m}\right).$$
Now we use the identity for the probabilists Hermite polynomials
$$ H e_n(x) =  \sum_{m=0}^{n/2} \frac{n!}{2^m m! (n-2m)!} x^{n-2m}$$
to rewrite the expression as
$$ S_1 = \frac{1}{n^{\ell/2}}\sum_{m=0}^{\ell} (-1)^m \binom{\ell}{m} He_m(\gamma) y^{\ell-m}.$$
The final ingredient is an addition formula for the (probabilists') Hermite polynomials
$$ He_n(a+b) = \sum_{k=0}^{n}{ \binom{n}{k} x^{n-k} He_k(y)}$$
allowing us to write
$$ S_1 =  \frac{(-1)^{\ell}}{n^{\ell/2}} He_{\ell}( \gamma-y),$$
and hence, using the symmetries of Hermite polynomials, we have
$$ {n^{\ell/2}}  S_1 = He_{\ell}( y-\gamma),$$
where $\gamma$ is a random variable. We recall that $\gamma \sim e_1/\sqrt{n}$ which implies that
$\gamma$ converges to a Gaussian distribution as $n \rightarrow \infty$. \end{proof}

\subsection{Proof of the Proposition}
\begin{proof}
Consider the polynomial $ (x^2-1)^n.$ We wish to characterize the behaviour of its first few derivatives in the vicinity of the origin. We begin by introducing the new variable $y = \sqrt{n}x,$ and observing that
$$(x^2-1)^n = (-1)^n\left(1-\frac{y^2}{n}\right)^n \sim (-1)^n e^{-y^2}.$$ 
This heuristic motivates the definition of the function $f_n,$ defined by
$$f_n(y)= e^{y^2}(1-y^2/n)^n.$$
 We observe that
$$f_n'(y) = 2ye^{y^2}(1-y^2/n)^n -2ye^{y^2}(1-y^2/n)^{n-1}=- \frac{1}{n}\frac{2y^3}{1-\frac{y^2}{n}}f_n(y).$$
We also observe that, for $|y| < \sqrt{n}$,
$$ \frac{2y^3}{1-\frac{y^2}{n}} = 2 y^3 \sum_{k=0}^{\infty} \frac{y^{2k}}{n^k} \qquad \mbox{as well as} \qquad e^{y^2} = \sum_{k=0}^{\infty} \frac{y^{2k}}{k!}$$
and
$$ \left(1 - \frac{y^2}{n}\right)^n = \sum_{k=0}^{n}{ (-1)^k \binom{n}{k} \frac{y^{2k}}{n^k} }.$$
We see that all three functions are power series with rapidly decaying coefficients. In particular, this implies that for any fixed $\ell \geq 1$ and any $a > 0$, there exists a constant $c_{a, \ell}$ such that for all $|x| \leq a$,
$$ |f_n^{(\ell)}(x)| \leq \frac{c_{\ell, a}}{n}.$$
Next, we let $\phi_{k,n}$ be the function defined by
$$\phi_{k,n} (y) = e^{-y^2} \frac{{\rm d}^k}{{\rm d}y^k}  f_n(y).$$
The above arguments guarantee that for all $|x| \leq a$, we have $|\phi_{k,n}(x)| \lesssim c_{a, \ell}^{*}/n$ for $k >0.$ Similarly, it follows from the definition of $f_n$ and $\phi_{0,n}$ that
$$\phi_{0,n}(y) = (1- y^2/n)^n.$$
Moreover, from the definition it is clear that
$$\phi_{k,n}'(x) = -2y \phi_{k,n}(y) + \phi_{k+1,n}(y).$$
A straightforward computation shows that
$$\phi_{0,n}^{(\ell)}(y) =  \sum_{j=0}^\ell (-1)^{\ell-j} {\ell \choose j} H_{\ell-j}(y) \phi_{j,n}(y),$$
which, together with the fact that $|\phi_{j,n}| \in O(1/n)$ for $j>0,$ shows that
$$\frac{{\rm d}^\ell}{{\rm d} y^\ell} \left(1 - \frac{y^2}{n} \right)^n = (-1)^\ell H_\ell(y) + O(1/n).$$
\end{proof}

\section{Comments and Remarks}

\subsection{The Moment Method} A classical approach to semicircle laws is the moment method: one computes all moments of the arising distribution and then deduces properties of the distribution from that.
We recall that, as derived above, we have
$$ \frac{\ell!}{n!} p_n^{(n-\ell)}(x) = \prod_{k=1}^{\ell}{(x - r_i)} =\sum_{m=0}^{\ell} (-1)^m f_m x^{\ell-m},$$
where 
$$ f_k = e_k  \frac{\ell!}{(\ell-k)!} \frac{(n-k)!}{n!}.$$ 
We recall that the $e_k$ are the elementary symmetric polynomials of the roots $r_i$. We will abbreviate their power sum as
$$ q_m = \sum_{i=1}^{\ell} r_i^m.$$
We observe that the power sums $q_m$ can be written in terms of the rescaled elementary symmetric polynomials via the following identity
$$ q_m = \sum_{s_1 + 2s_2 + \dots + m s_m = m \atop s_1 \geq 0, s_2 \geq 0, \dots, s_m \geq 0} (-1)^m \frac{m(s_1 + s_2 + \dots + s_m - 1)!}{s_1! s_2! \dots s_m!} \prod_{i=1}^{m}(-f_i)^{s_i}.$$
This seems like it would lead to some combinatorial identities: we expect the roots $r_i$ to follow a semicircle distribution as $\ell$ becomes large allowing us to approximate $q_m$ by the $\ell-$th moment of a semicircle distribution (suitably scaled). Conversely, using the Lemma proved above, we can approximate $f_k$ to leading order by a suitably-scaled Hermite polynomial. This seems reminiscent of work of Carlitz \cite{carl}.

\subsection{A Connection to Random Projections} The behaviour of the roots of polynomials after differentiation has a natural and classical connection to changes in eigenvalues after restriction to certain codimension one subspaces. In particular, if $A$ is an $n \times n$ diagonal matrix with entries $\lambda_1,\dots,\lambda_n$ and $P$ is the projection matrix $\frac{1}{n} \mathbf{1} \mathbf{1}^T,$ where $\mathbf{1}$ is the vector of all ones, then the non-zero eigenvalues of $(I-P)A(I-P)$ are the solutions of the following equation
$$0 = \sum_{i=1}^n \frac{1}{z-\lambda_i},$$
and correspond to the roots of the derivative of the characteristic polynomial of $A.$ 

   \begin{figure}[h!]
    \centering
        \includegraphics[width=0.7\textwidth]{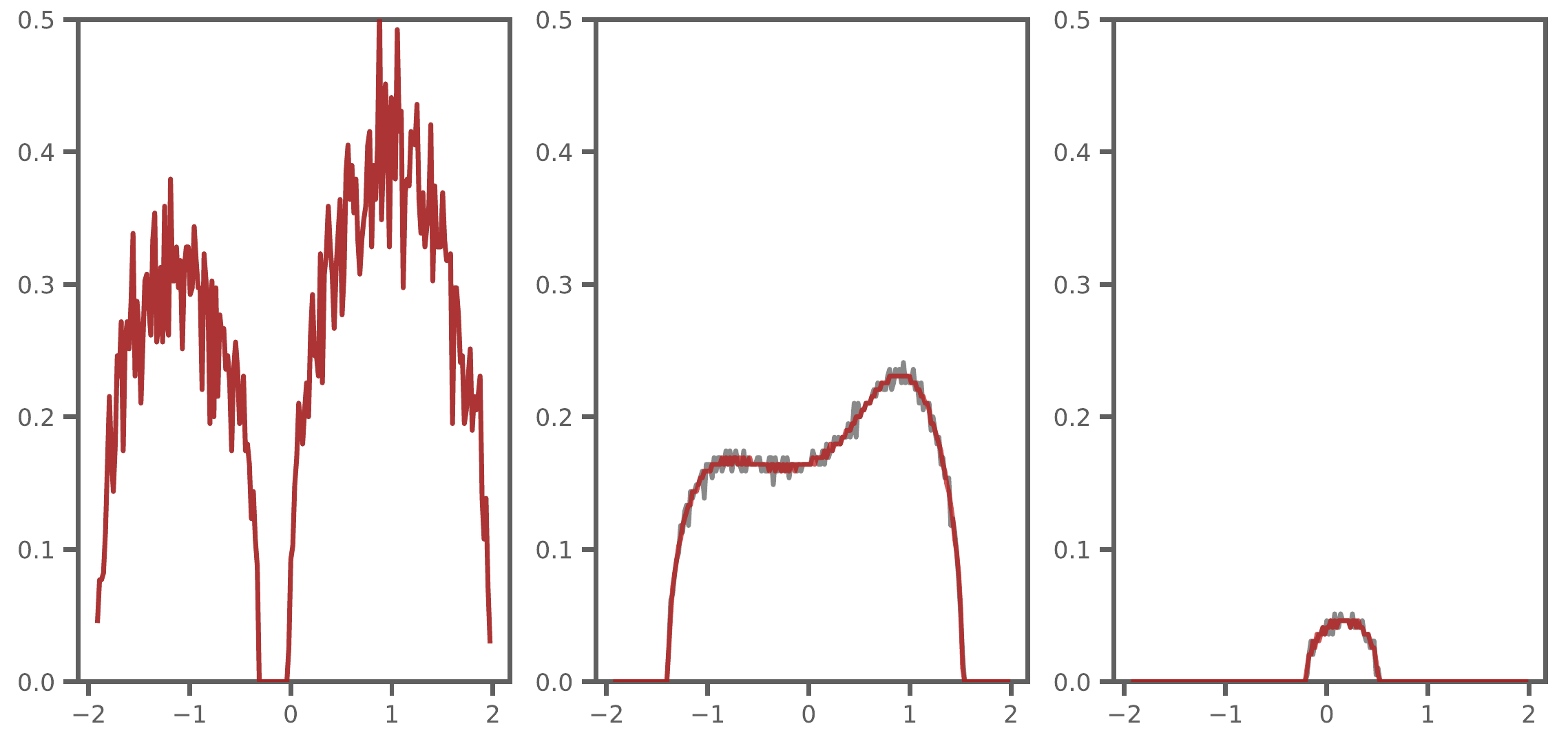}
    \caption{Histograms of eigenvalues after successive projections, obtained using the deterministic (red) and black (random) eigenvalue update formulae. left: the initial (1000) eigenvalues; center: the eigenvalues after 500 projections; right: the eigenvalues after 975 projections. }\label{fig:projs}
\end{figure}

An interesting and related question is what happens if a random projection is used instead? Namely, if $A$ is an $n \times n$ real symmetric matrix, $\|v\| = 1$ is uniformly distributed on $S^{n-1}$ and $P_v = vv^T$, then what are the eigenvalues of $(I-P_v) A (I-P_v)?$ Applying the Sherman-Morrison-Woodbury formula \cite{bunch,golub,gu,simonbook,sherman,woodbury} shows that the new eigenvalues are the roots of the equation
$$0 = \sum_{i=1}^n \frac{|w_i|^2}{z-\lambda_i},$$
where $(w_1,\dots,w_n)$ is uniformly distributed on $\mathbb{S}^{n-1}.$ Clearly, one might expect that the roots of this equation will roughly coincide with those of the deterministic case. It is not quite as obvious what will happen after this process of projection is iterated. Numerical experiments (as in Figure \ref{fig:projs}) indicate that the dynamics are related. If this is true, then the main result of this paper has immediate implications for the restriction of symmetric matrices to low-dimensional random subspaces.

\subsection{A concluding example.}

  \begin{figure}[h!]
  
    \centering
    \begin{subfigure}[b]{\textwidth}
        \includegraphics[width=0.9\textwidth]{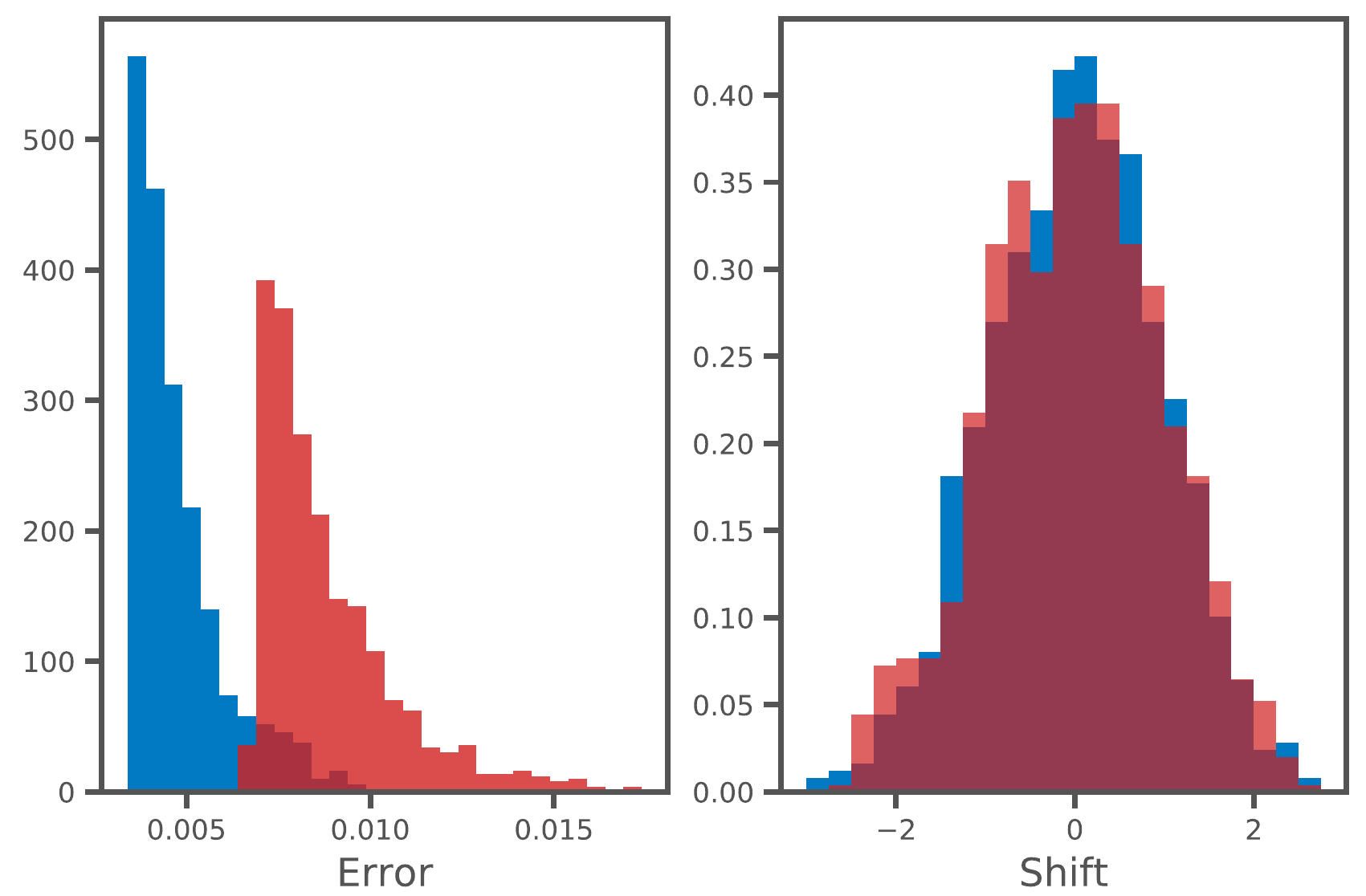}
        \label{fig:flat}
    \end{subfigure}\\
        ~ 
    \begin{subfigure}[b]{0.4\textwidth}
        \includegraphics[width=\textwidth]{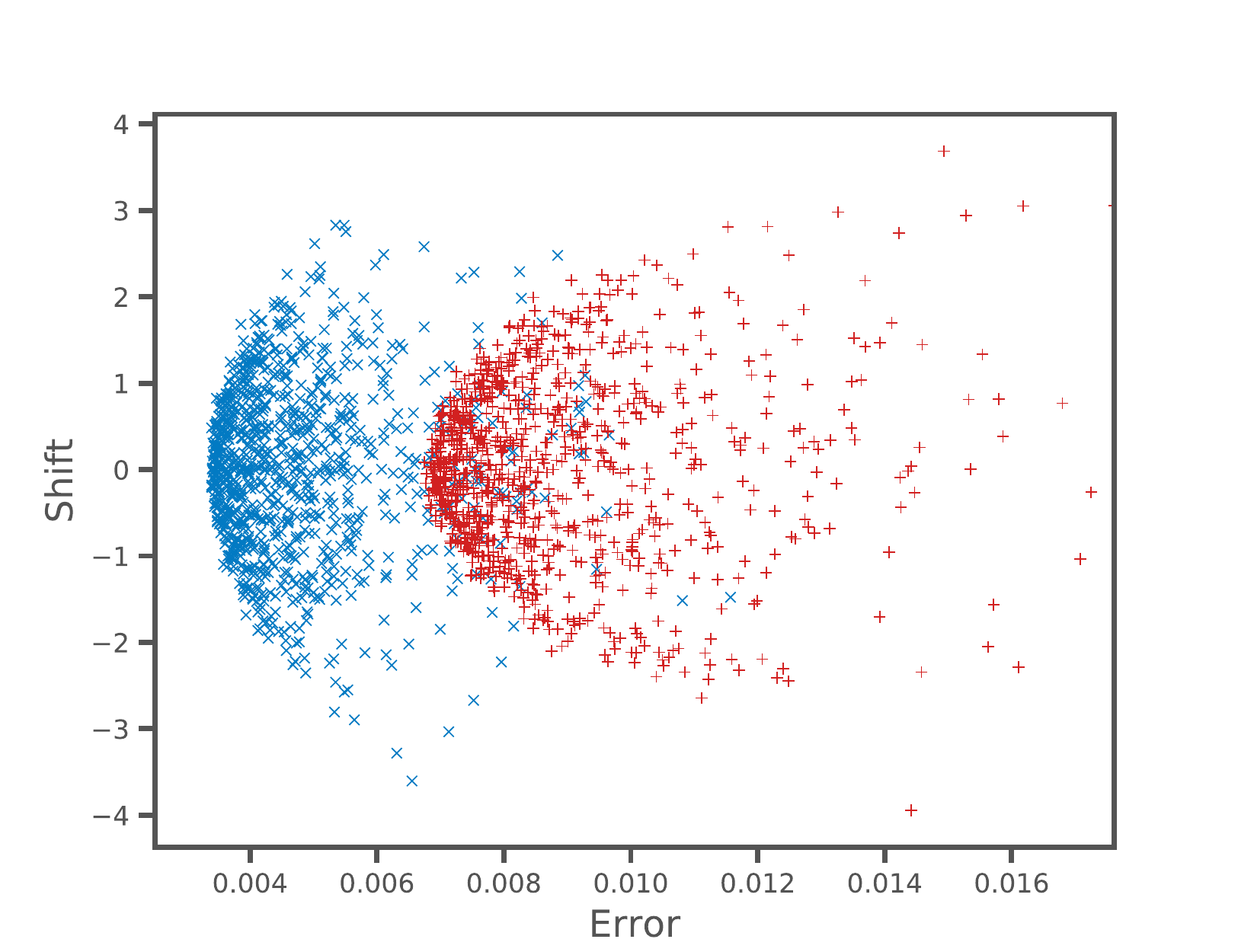}
        \label{fig:parab}
    \end{subfigure}
        \begin{subfigure}[b]{0.49\textwidth}
        \includegraphics[width=\textwidth]{fort_polys}
        \label{fig:gauss}
    \end{subfigure}
\captionsetup{singlelinecheck=off}
\end{figure}

We conclude with a more detailed numerical example given by the initial probability density
\[
   \mu(x) = \frac{9 \sqrt{3} } {10\sqrt{5}} x^2\,\chi_{[-\sqrt{5/3},\sqrt{5/3}]},
\]
We compute examples comprised of 1000 initial roots (in red) and 2000 initial roots (in blue). The first picture (upper left) in Figure 7 shows
histograms of the error $$\frac{1}{50}\left(\sum_{i=1}^{50} (\sqrt{n}\, r_i+\gamma- y_i)^2\right)^\frac{1}{2}$$ where $r_i$ are the roots of $p_n^{(n-50)},$ $y_i$ are the roots of $He_{50},$ and $\gamma$ is the best shift. The second picture (upper right) shows histograms of the shifts $\gamma$ selected to match $\sqrt{n}\,  r_i$ with $y_i$. They can be seen to be close to a Gaussian. The third picture (lower left) shows a scatter plot of the error and the shifts and the last figure (lower right) shows a plot of 
$$\log_{10} \left| 1000^{25} \frac{50!}{1000!} \, p_{1000}^{950}\left( \frac{x}{\sqrt{n}} \right)\right| \quad \mbox{(in purple)}$$
 compared to $\log_{10} \left| He_{50}(x)\right|$ (in green).

\newpage

\end{document}